\documentclass[mathpazo]{cicp}
\usepackage{ifpdf}
\ifpdf \else \pdfoutput=1 \fi

\usepackage{graphicx, epsfig, psfrag}
\usepackage{amsfonts,amsmath,amssymb,amsbsy,amsthm}
\usepackage{framed, multirow}
\usepackage{latexsym}
\usepackage{url}
\usepackage{xcolor}
  \definecolor{newcolor}{rgb}{.8,.349,.1}
  \definecolor{b+}{rgb}{0.5, 0.5, 1}
  \definecolor{b-}{rgb}{0, 0, 0.5}
  \definecolor{r|b-}{rgb}{0.5, 0, 0.25}
  \definecolor{g|b-}{rgb}{0, 0.5, 0.25}
  \definecolor{shadecolor}{rgb}{0.75, 0.75, 1}
\usepackage[
  colorlinks,
  urlcolor = black,
  linkcolor = b-,
  anchorcolor = r|b-,
  citecolor = g|b-,
  hypertexnames = false % Just add this command to avoid warning like `destination with the same identifier has been already used, duplicate ignored'. The problem may come from the `float' package.
]{hyperref}
\usepackage{geometry}
\usepackage{amsthm, amsmath, amssymb, amsfonts, amsbsy}
\usepackage{mathrsfs} % to use \mathscr
\usepackage{mathtools}
\allowdisplaybreaks
\usepackage{indentfirst}
\usepackage{graphicx}
\usepackage{epstopdf} % epstopdf after graphicx
\usepackage{makeidx}
  \makeindex

\usepackage{pdfpages} % using \includepdf [addtotoc={1,section,1,title in toc,cc},pages=1-7,offset=0cm 0.5cm] {mypdffile.pdf}
\usepackage{subfigure}
\usepackage{longtable}
\usepackage{hyperref}
\usepackage{esint}
\usepackage{soul}
\usepackage{float}
\PassOptionsToPackage{hyphens}{url}\usepackage{hyperref}

\def\XXint#1#2#3{{\setbox0=\hbox{$#1{#2#3}{\int}$ }
\vcenter{\hbox{$#2#3$ }}\kern-.6\wd0}}
\def\<{\langle}
\def\>{\rangle}

\begin{document}
%%%%% title : short title may not be used but TITLE is required.
% \title{TITLE}
% \title[short title]{TITLE}
\title{Blended Ghost Force Correction Method for 3D Crystalline Defects}

%%%%% author(s) :
% single author:
% \author[name in running head]{AUTHOR\corrauth}
% [name in running head] is NOT OPTIONAL, it is a MUST.
% Use \corrauth to indicate the corresponding author.
% Use \email to provide email address of author.
% \footnote and \thanks are not used in the heading section.
% Another acknowlegments/support of grants, state in Acknowledgments section
% \section*{Acknowledgments}
% \author[O.~Author]{Only Author\corrauth}
% \address{School of Mathematical Sciences, Beijing Normal University,
% Beijing 100875, P.R. China}
% \email{{\tt author@email} (O.~Author)}

\author[L. Fang, L. Zhang]{Lidong Fang\affil{1}\comma, and Lei Zhang\affil{1}\corrauth}
\address{\affilnum{1}\ School of Mathematical Sciences, Institute of Natural Sciences, and MOE-LSC, Shanghai Jiao Tong University, Shanghai, 200240, China.}
\emails{{\tt ldfang.sjtu@gmail.com} (L. Fang), {\tt lzhang2012@sjtu.edu.cn} (L. Zhang)}

%%%%% Begin Abstract %%%%%%%%%%%
\begin{abstract}
Atomistic/continuum coupling method is a class of multiscale computational method for the efficient simulation of crystalline defects. The recently developed blended ghost force correction (BGFC) method combines the efficiency of blending methods and the accuracy of QNL type methods. BGFC method can be applied to multi-body interaction potentials and general interfaces. In this paper, we present the formulation, implementation and analysis of the BGFC method in three dimensions. In particular, we focus on the difference and connection with other blending variants, such as energy based blended quasi-continuum method (BQCE) and force based blended quasi-continuum method (BQCF). The theoretical results are justified by a few benchmark numerical experiments with point defects and microcrack in the three dimensional FCC lattice. 
\end{abstract}
%%%%% end %%%%%%%%%%%

%%%%% AMS/PACs/Keywords %%%%%%%%%%%
%\pac{}
%\ams{52B10, 65D18, 68U05, 68U07}

%%%% maketitle %%%%%
\maketitle

\newpage



\section*{Supplementary material}

# 1  Introduction

Atomistic/continuum (a/c) coupling method is a class of computational multiscale methods [36] that aim to combine the accuracy of fine scale models and the efficiency of coarse scale models for crystalline defects. Namely, fine scale models can be applied in a small neighborhood of the localized defects such as vacancies, dislocations, and cracks, while coarse scale models can be employed away from the defect cores where elastic deformation occurs.

In the past two decades, a/c methods have attracted great attention from both the engineering community and the mathematical community [1, 2, 13, 18, 22, 24, 32, 35, 36]. On one hand, predictive simulations for materials defects such as point defects and dislocations are essential to underpin the elastic and plastic deformation mechanism of materials [28]; on the other hand, the quantitative estimates of the approximation error for a/c methods as a representative concurrent multiscale method help elucidate open questions and establish an analytical framework for similar multiscale computational methods [3, 7, 12, 29].

For a/c coupling methods, fine scale models are usually empirical interaction potentials, while coarse scale models are coarse-grained continuum elastic models. Energy based methods and force based methods are two major classes of a/c coupling methods, we refer to [18, 22, 36] for reviews of many existing a/c methods. Energy based methods construct a hybrid energy functional as a weighted combination of atomistic and continuum energy functional, and one of the major challenges for energy based methods is to eliminate the so-called "ghost forces" [22] near the atomistic/continuum interface. Force based methods compute the equilibrium of atomistic and continuum forces from corresponding energies, see [6,16,17,21,23] for recent advances. In practice, force based methods can remove ghost forces, and seem to be optimal in terms of coupling error. However, they are not conservative, namely, there is no associated energy functional, and it is usually difficult to establish the stability of the force operator.

Blending type energy based a/c coupling methods [15, 19, 37] smear out the a/c interface and thus propose a weighted energy functional by a blending function. It is easy to implement, however, it does not eliminate the ghost forces as consistent methods do, and only has suboptimal convergence rate [15].

QNL (quasi-nonlocal) type energy based a/c coupling methods aim to eliminate the ghost force. Therefore, they are referred to as consistent coupling methods. QNL type methods have been well developed in 1D and 2D for multi-body interactions and general interfaces [8, 26, 33]. However, in three dimensions, QNL type methods are only available for pair interactions [20, 30, 31], and the construction for multi-body interactions remains open.

In [27], we constructed the blended ghost force correction (BGFC) scheme by integrating two popular ideas: blending [37] and ghost force correction [32]. The BGFC scheme combines the efficiency of the blending methods as well as the accuracy of the QNL type consistent methods. It is quasi-optimal in the sense that it yields the same convergence rate as the force based a/c coupling schemes [14–16]. In fact, it is most

instructive to derive the scheme through a modification of the site energies, which can be regarded as a prediction-correction scheme.

For simplicity, the implementation of the BGFC method in [27] is restricted in two dimensions. In this paper, we extend the BGFC scheme to three dimensions. This extension is highly nontrivial since the lattice structure, the partition of the graded mesh, and the implementation of the finite element method are much more complicated in three dimensions. We implement the 3D BGFC method for the FCC lattice and a second-nearest-neighbor multi-body interaction potential, and test a few prototypical benchmark examples such as the single-vacancy, separated-vacancy, and microcrack.

The paper is organized as follows. In Section 2, we briefly describe the model setup, formulate the BGFC method, and also state the main theoretical results. In Section 3, we discuss the numerical implementation in 3D, and demonstrate the convergence of the BGFC scheme with typical numerical examples in the three dimensional FCC lattice. In Section 4, we sketch the analytical proof for the 3D BGFC method, which focuses on the difference of BGFC and other blending variants such as energy based blended quasi-continuum (BQCE) and force based blended quasi-continuum (BQCF). We conclude the paper and discuss future research directions in Section 5 .

**Notation**

We use the symbol $\langle \cdot, \cdot \rangle$ to denote an abstract duality pair between a Banach space $V$ and its dual space $V^*$. For second order tensors $A$ and $B$, we denote $A:B=\sum_{i,j}A_{ij}B_{ij}$ and $A \otimes B$ the standard kronecker product. For functional $E \in C^2(X)$, the first and second variations are denoted by $\langle \delta E(u), v \rangle$ and $\langle \delta^2 E(u)v, w \rangle$ for $u, v, w \in X$, respectively. For a finite set $A$, we will use $\#A$ to denote the cardinality of $A$. The closed ball with radius $r$ and center $x$ is denoted by $B_r(x)$, or $B_r$ if the center is the origin.

## 2 BGFC Method

### 2.1 Atomistic lattice and lattice functions

Given $d \in \{2,3\}$, $\mathsf{A} \in \mathbb{R}^{d \times d}$ non-singular, $\Lambda^{\text{hom}} := \mathsf{A}\mathbb{Z}^d$ is the homogeneous reference lattice which represents a perfect single lattice crystal formed by identical atoms. For face centered crystal (FCC) and body centered crystal (BCC) lattices, we have

$$\mathsf{A}_{\text{FCC}} = \frac{1}{2} \begin{pmatrix} 0 & 1 & 1 \\ 1 & 0 & 1 \\ 1 & 1 & 0 \end{pmatrix}, \quad \mathsf{A}_{\text{BCC}} = \frac{1}{2} \begin{pmatrix} -1 & 1 & 1 \\ 1 & -1 & 1 \\ 1 & 1 & -1 \end{pmatrix}, \tag{2.1}$$

and we note that $\mathsf{A}$ is not unique.

$\Lambda \subset \mathbb{R}^d$ is the reference lattice with some local defects. The mismatch between $\Lambda$ and $\Lambda^{\text{hom}}$ represents possible defects $\Lambda^{\text{def}}$, which are contained in some localized defect cores $\Omega^{\text{def}}$ with radius $R^{\text{def}} > 0$ such that the atoms in $\Lambda \setminus \Omega^{\text{def}}$ do not interact with defects $\Lambda^{\text{def}}$. Vacancy, interstitial and impurity are different types of possible point defects.

### 2.1.1 Lattice function and lattice function space

Denote the set of vector-valued lattice functions by

$$\mathscr{U} := \{v : \Lambda \to \mathbb{R}^d\}.$$

and $\mathscr{U}^c := \{u | \text{supp}(u) \text{ is compact}\}$.

A deformed configuration is a lattice function $y \in \mathscr{U}$. Let $x$ be the identity map, the displacement $u \in \mathscr{U}$ is defined by $u(\ell) = y(\ell) - x(\ell) = y(\ell) - \ell$ for any $\ell \in \Lambda$.

We can introduce the discrete homogeneous Sobolev spaces

$$\mathscr{U}^{1,2} := \{v : \Lambda \to \mathbb{R}^d \,|\, |v|_{\mathscr{U}} < +\infty\},$$

with semi-norm $|v|_{\mathscr{U}^{1,2}} := \left(\sum_{\ell \in \Lambda} \sum_{\ell' \in \mathcal{N}_\ell} |v(\ell') - v(\ell)|^2\right)^{1/2}$, where $\mathcal{N}_\ell$ is the set of nearest neighbors of $\ell$ in $\Lambda$.

## 2.2 Atomistic Model and Continuum Model

### 2.2.1 Atomistic Model

For $\ell \in \Lambda^{\text{hom}}$, the interaction range $\mathcal{R}_\ell := \{\ell' - \ell | 0 < |\ell' - \ell| \leq R^{\text{cut}}, \ell' \in \Lambda^{\text{hom}}\}$ is the union of lattice vectors defined by the finite differences within radius $R^{\text{cut}}$. We define the "finite difference stencil" as $Dv(\ell) := \{D_\rho v(\ell)\}_{\rho \in \mathcal{R}_\ell} := \{v(\ell+\rho) - v(\ell)\}_{\rho \in \mathcal{R}_\ell}$.

Let $V_\ell(Du)$ denote the site energy associated with the lattice site $\ell \in \Lambda^{\text{hom}}$. We assume that the potential $V_\ell(Du) \in C^k((\mathbb{R}^d)^{\mathcal{R}_\ell})$, $k \geq 2$. We also assume that $V_\ell(Du)$ is homogeneous outside the defect region $\Omega^{\text{def}}$, namely, $V_\ell = V$ and $\mathcal{R}_\ell = \mathcal{R}$ for $\ell \in \Lambda \setminus \Omega^{\text{def}}$. Furthermore, $V$ and $\mathcal{R}$ have the following point symmetry: $\mathcal{R} = -\mathcal{R}$, and $V(\{-g_{-\rho}\}_{\rho \in \mathcal{R}}) = V(g)$.

In the following, for clarity of presentation, we mainly consider the single vacancy case and assume that the vacancy locates at the origin, namely, $\Lambda^{\text{hom}} \setminus \Lambda = \{0\}$. The formulation can be generalized to separated defects, microcrack, and dislocations [27]. We will be show the numerical experiments for single vacancy, separated vacancies, and microcrack in Section 3.

The energy of an infinite configuration is typically ill-defined. We can introduce the *renormalized potential* as $V'_\ell(Du) := V_\ell(Du) - V_\ell(Du_0)$ for some reference displacement $u_0 \in \mathscr{U}^{1,2}$ (for point defects, we can simply choose $u_0 = 0$), and define the energy-difference as

$$\mathscr{E}^{\text{a}}(u) = \sum_{\ell \in \Lambda^{\text{hom}}} V'_\ell(Du), \tag{2.2}$$

and the *defect potential* as

$$\mathscr{P}^{\text{def}}(u) = -\sum_{\ell \in \Lambda^{\text{hom}} \setminus \Lambda} V'_\ell(Du), \tag{2.3}$$

where we can assume $u(\ell) = 0$, for $\ell \in \Lambda^{\rm hom} \setminus \Lambda$. We can take the site potential for $\ell \in \Lambda^{\rm hom} \setminus \Lambda$ as $V_\ell \equiv 0$ for vacancies, and as certain nonzero potential for interstitials or impurities. We note that $\mathscr{P}^{\rm def}$ is localized ($\mathscr{P}^{\rm def}(u)$ only depends on $(u(\ell), |\ell| \leq R^{\rm def})$) and translation invariant ($\mathscr{P}^{\rm def}(u) = \mathscr{P}^{\rm def}(u+c)$, for a constant $c$).

The atomistic problem now reads

$$u^{\rm a} \in \arg\min\{\mathscr{E}^{\rm a}(v) + \mathscr{P}^{\rm def}(v) \,|\, v \in \mathscr{U}^{1,2}\}, \tag{2.4}$$

where "argmin" is understood as the set of local minima.

We call $u^{\rm a}$ a strongly stable solution to (2.4) if there exists $\gamma_{\rm a} > 0$ such that

$$\langle \delta[\mathscr{E}^{\rm a} + \mathscr{P}^{\rm def}](u^{\rm a}), v \rangle = 0 \quad \text{and} \quad \langle \delta^2[\mathscr{E}^{\rm a} + \mathscr{P}^{\rm def}](u^{\rm a})v, v \rangle \geq \gamma_{\rm a} |v|_{\mathscr{U}^{1,2}}^2 \quad \forall v \in \mathscr{U}^{1,2}. \tag{2.5}$$

It is shown that the energy-difference functional $\mathscr{E}^{\rm a} + \mathscr{P}^{\rm def}$ is well-defined under suitable conditions (regularity and homogeneity outside the defect) [10, Lemma 2.1].

### 2.2.2 Continuum model

A continuum model can be derived by coarse graining the corresponding atomistic model, and computationally it allows for the reduction of degrees of freedom when the deformation is smooth. Cauchy-Born continuum model is a typical choice in the multi-scale context [9, 25]. The Cauchy-Born energy density $W: \mathbb{R}^{d \times d} \to \mathbb{R}$ is defined by

$$W(\mathsf{F}) := \det(\mathsf{A}^{-1}) V(\mathsf{F}\mathcal{R}).$$

and the Cauchy-Born energy difference is defined by

$$W'(\mathsf{F}) := W(\mathsf{F} + \mathsf{I}) - W(\mathsf{I}), \quad \forall \mathsf{F} \in \mathbb{R}^{d \times d}.$$

## 2.3 BGFC model

In this section, we introduce the BGFC model, as well as related blended coupling variants such as BQCE and BQCF models.

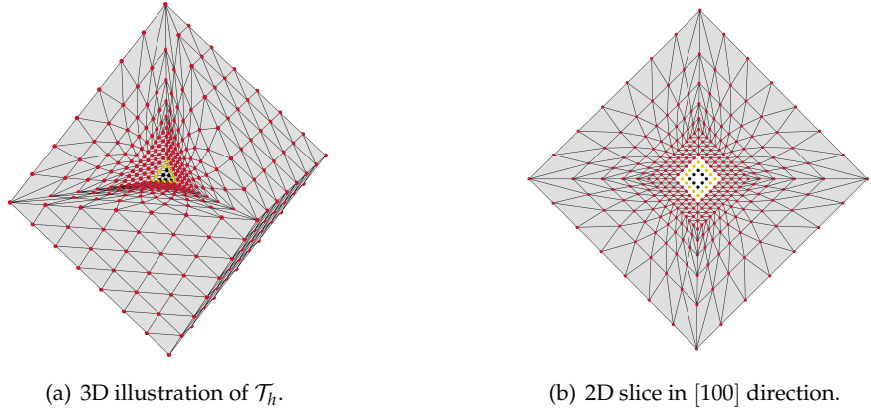

(a) 3D illustration of $\mathcal{T}_h$.   (b) 2D slice in [100] direction.

Figure 1: Example of the coarse grained mesh $\mathcal{T}_h$ for a point vacancy in the 3D FCC lattice. The nodes with $\beta=1$ are colored red, the nodes with $\beta=0$ are colored black, and the nodes in the blending region are colored yellow.

To construct the blended schemes, we define a regular simplicial finite element grid $\mathcal{T}_h$ with nodes $\mathcal{X}_h$, with the minimal requirement that $\mathcal{X}_h \cap B_{R^{\text{def}}} = \Lambda \cap B_{R^{\text{def}}}$, namely, the defect core is resolved exactly. See Figure 1 for an illustration of $\mathcal{T}_h$ for a point vacancy in 3D FCC lattice. Let $\text{DoF} := \#\mathcal{X}_h$, and $\Omega_h := \bigcup \mathcal{T}_h$ be the resulting computational domain. There exists constants $0 < R^{\text{i}} < R^{\text{o}}$, such that $B_{R^{\text{i}}} \subset \Omega_h \subset B_{R^{\text{o}}}$.

$\Omega_h$ can be partitioned into the atomistic region $\Omega^{\text{a}}$ and the continuum region $\Omega^{\text{c}}$, outwards from the defect. We introduce the blending function $\beta \in C^{2,1}(\mathbb{R}^d)$ with $\beta = 0$ in $B_{R^{\text{a}}} \subset \Omega^{\text{a}}$ and $\beta = 1$ in $\mathbb{R}^d \setminus \Omega^{\text{a}}$, where $R^{\text{def}} < R^{\text{a}} < R^{\text{i}}$. The blending region $\Omega^{\text{b}} := \text{supp}(\nabla \beta) + B_{2R^{\text{cut}} + \sqrt{d}} = \Omega^{\text{a}} \cap \Omega^{\text{c}}$ has width $R^{\text{b}} - R^{\text{a}} \simeq R^{\text{a}}$. The exterior region is $\Omega^{\text{ext}} := \mathbb{R}^d \setminus B_{R^{\text{i}}/2}$. See Figure 2 for an illustration of the partition in 2D. $(\beta, \mathcal{T}_h)$ are the main approximation parameters for the blended coupling methods, and can be characterized by the above parameters $(R^{\text{a}}, R^{\text{b}}, R^{\text{i}}, R^{\text{o}}, \beta, h(x), \text{DoF})$.

We assume that in the continuum region $\Omega^{\text{c}}$, the finite element tetrahedra mesh is graded with a mesh size function $h(x) = \text{diam}(T)$ for $x \in T \subset \mathcal{T}_h$. Let the space of coarse-grained admissible displacements be given by

$$\mathscr{U}_h := \{ v_h \in C(\mathbb{R}^d; \mathbb{R}^d) \,|\, v_h \text{ is p.w. affine w.r.t. } \mathcal{T}_h, \text{ and } v_h|_{\mathbb{R}^d \setminus \Omega_h} = 0 \}.$$

Let $Q_h$ denote the $P_0$ midpoint interpolation operator, so that $\int_{\Omega_h} Q_h f$ is the midpoint rule approximation to $\int_{\Omega_h} f$. We first formulate the energy based blended quasicontinuum (BQCE) scheme as introduced in [19] and analyzed in [15]. Let $V'_\ell$ be the renormalized potential defined in (2.2) and assume homogeneity $V'_\ell = V'$ outside the defect core, we define the BQCE energy functional for $u_h \in \mathscr{U}_h$ as

$$\mathscr{E}_h^{\text{bqce}}(u_h) := \sum_{\ell \in \Lambda^{\text{hom}} \cap \Omega_h} (1 - \beta(\ell)) V'_\ell(Du) + \int_{\Omega_h} Q_h \big[ \beta W'_\ell(\nabla u_h) \big]. \tag{2.6}$$

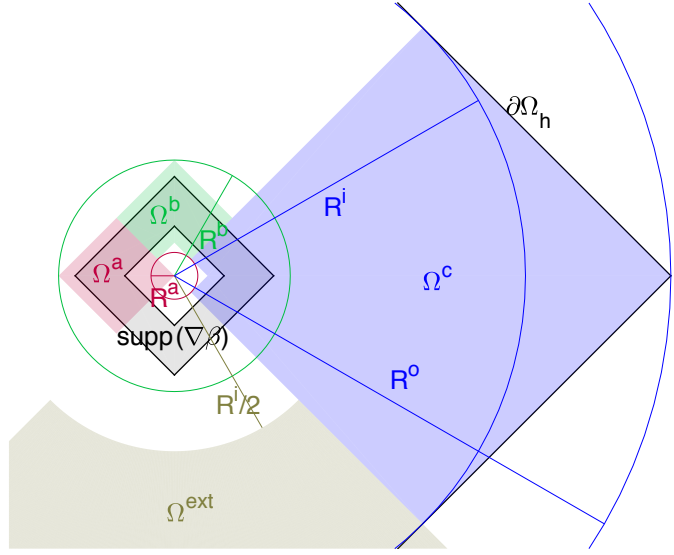

Figure 2: 2D illustration for the partition of the computational domain $\Omega_h$, which shows the domains $\Omega^a$, $\Omega^b$, $\Omega^c$, $\Omega^{\text{ext}}$, and radii of domains $R^a$, $R^b$, $R^i$, $R^o$.

The BQCE problem is to compute

$$u_h^{\text{bqce}} \in \arg\min_{u_h \in \mathscr{U}_h} \left\{ \mathscr{E}_h^{\text{bqce}}(u_h) + \mathscr{P}^{\text{def}}(u_h) \right\}, \tag{2.7}$$

where $u_h^{\text{bqce}}$ is a BQCE solution.

The BGFC formulation is based on a *second renormalization* of the potential. which can be defined by

$$V_\ell''(Du) := V_\ell(Du) - V_\ell(0) - \langle \delta V_\ell(0), Du \rangle,$$

for all $u \in \mathscr{U}^{1,2}$. We take the reference solution as $u_0 = 0$ for the single vacancy case.

The corresponding second renormalized Cauchy-Born energy density is

$$W''(\mathsf{F}) := W(\mathsf{F}+\mathsf{I}) - W(\mathsf{I}) - \langle \delta W(\mathsf{I}), \mathsf{F} \rangle,$$

for $\mathsf{F} \in \mathbb{R}^{d \times d}$. It follows that $\partial V_\ell''(0) = 0$ and $\partial W''(0) = 0$.

The BGFC energy for $u_h \in \mathscr{U}_h$ is defined by

$$\mathscr{E}_h^{\text{bgfc}}(u_h) := \sum_{\ell \in \Lambda^{\text{hom}} \cap \Omega_h} (1 - \beta(\ell)) V_\ell''(u_h) + \int_{\Omega_h} Q_h \big[ \beta W''(\nabla u_h) \big]. \tag{2.8}$$

The BGFC problem is to compute

$$u_h^{\text{bgfc}} \in \arg\min_{u_h \in \mathscr{U}_h} \left\{ \mathscr{E}_h^{\text{bgfc}}(u_h) + \mathscr{P}^{\text{def}}(u_h) \right\}, \tag{2.9}$$

where $u_h^{\text{bgfc}}$ is a BGFC solution. We denote the BGFC scheme with piecewise linear space $\mathscr{U}_h$ as P1-BGFC.

We can apply P2 finite elements in the continuum region to achieve optimal convergence rate, see Theorem 2.1 and the discussions in Section 2.4. We decompose $\mathcal{T}_h = \mathcal{T}_h^{(\text{P1})} \cup \mathcal{T}_h^{(\text{P2})}$, where $\mathcal{T}_h^{(\text{P1})} = \{T \in \mathcal{T}_h \,|\, T \cap \Omega^{\text{a}} \neq \emptyset\}$, and define

$$\mathscr{U}_h^{(2)} := \big\{u_h \in C(\mathbb{R}^d;\mathbb{R}^m) \,\big|\, u_h|_T \text{ is affine for } T \in \mathcal{T}_h^{(\text{P1})},$$
$$u_h|_T \text{ is quadratic for } T \in \mathcal{T}_h^{(\text{P2})}, \text{ and} \qquad (2.10)$$
$$u_h = 0 \text{ in } \mathbb{R}^d \setminus \Omega_h\big\}.$$

By adjusting the quadrature operator $Q_h$ such that $\nabla u_h \otimes \nabla u_h$ can be integrated exactly for $u_h \in \mathscr{U}_h^{(2)}$, we define the P2-BGFC problem by

$$u_h^{\text{bgfc},2} \in \arg\min_{u_h \in \mathscr{U}_h^{(2)}} \big\{\mathscr{E}_h^{\text{bgfc}}(u_h) + \mathscr{P}^{\text{def}}(u_h)\big\}. \qquad (2.11)$$

We also introduce the BQCF scheme here. The BQCF operator is the nonlinear map $\mathscr{F}_h^{\beta}: \mathscr{U}_h \to \mathscr{U}_h^*$, defined by

$$\big\langle \mathscr{F}_h^{\beta}(u_h), v_h \big\rangle := \big\langle \delta\mathscr{E}^{\text{a}}(u_h), (1-\beta)v_h \big\rangle + \big\langle \delta\mathscr{E}^{\text{c}}(u_h), \overline{\beta v_h} \big\rangle, \qquad (2.12)$$

where $\overline{\beta v_h}$ is the P1 interpolant of $\beta v_h$ (see Section 6.1), and $\mathscr{E}^{\text{c}}(u_h) := \int_{\Omega_h} Q_h[W(\nabla u_h + \mathsf{I}) - W(\mathsf{I})]$.

For the BQCF method, we approximate the atomistic problem (2.4) by the variational nonlinear system for $u_h^{\text{bqcf}} \in \mathscr{U}_h$:

$$\big\langle \mathscr{F}_h^{\beta}(u_h^{\text{bqcf}}) + \mathscr{P}^{\text{def}}(u_h^{\text{bqcf}}), v_h \big\rangle = 0, \quad \forall v_h \in \mathscr{U}_h. \qquad (2.13)$$

We also refer to BQCF with $\mathscr{U}_h$ as P1-BQCF. If one replaces $\mathscr{U}_h$ by $\mathscr{U}_h^{(2)}$ in (2.13), we obtain the P2-BQCF method.

### 2.4 Convergence Theorem

To measure the local "regularity" of a displacement function $u \in \mathscr{U}^{1,2}$, we can introduce a piecewise linear interpolation $\bar{u}$, and a $C^{2,1}$-conforming interpolant $\tilde{u}$ with respect to the atomistic grid $\Lambda^{\text{hom}}$, see Section 6.1 and Section 6.2 for more details. With the help of those interpolants, we are able to define the following error contributions in the convergence theorem Theorem 2.1.

- $E^{\text{apx}}$, best-approximation error:

$$E^{\text{apx}}(u) := \|\nabla \bar{u}\|_{L^2(\Omega^{\text{ext}})} + \|h\nabla^2 \tilde{u}\|_{L^2(\Omega^{\text{c}})} + \|h^2 \nabla^3 \tilde{u}\|_{L^2(\Omega^{\text{c}})}. \qquad (2.14)$$

$$E^{\text{apx},2}(u):=\|\nabla\bar{u}\|_{L^2(\Omega^{\text{ext}})}+\|h^2\nabla^3\tilde{u}\|_{L^2(\Omega^c)}. \tag{2.15}$$

where the first term $\|\nabla\bar{u}\|_{L^2(\Omega^{\text{ext}})}$ is the far field truncation error, and the remaining terms represent the finite element discretization/coarsening error.

- $E^{\text{cb}}$, Cauchy-Born modeling error:

$$E^{\text{cb}}(u):=\|\nabla^3\tilde{u}\|_{L^2(\Omega^c)}+\|\nabla^2\tilde{u}\|_{L^4(\Omega^c)}^2. \tag{2.16}$$

- $E^{\text{int}}$, coupling/interface error:

$$E^{\text{int}}(u):=\|\nabla^2\beta\|_{L^2(\Omega^b)}+\|\nabla\beta\|_{L^\infty(\Omega^b)}\|\nabla^2\tilde{u}\|_{L^2(\Omega^b)}, \tag{2.17}$$

$$E^{\text{int},2}(u):=\|\nabla^2\beta\|_{L^2(\Omega^b)}\|\nabla\tilde{u}\|_{L^2(\Omega^b)}+\|\nabla\beta\|_{L^\infty(\Omega^b)}\|\nabla^2\tilde{u}\|_{L^2(\Omega^b)}. \tag{2.18}$$

The (optimal) choice for the approximation parameters in $(\beta,\mathcal{T}_h)$ are given in the following Assumption 1.

**Assumption 1.** There exist regularity constants $\mathbf{C}=(C_b,C_\Omega,C_h)$ for $(\beta,\mathcal{T}_h)$, such that

$$R^b\leq C_b R^a, \quad \|\nabla^j\beta\|_{L^\infty}\leq C_b(R^a)^{-j}, j=1,2,3,$$
$$h(x)=\mathcal{O}(\max\{1,(|x|/R^a)^\chi\}), \quad 1<\chi<1+d/4,$$
$$\text{DoF}\leq C_h(R^a)^d\log(R^a), \quad \max_{T\in\mathcal{T}_h}h_T^d/|T|\leq C_h, \tag{2.19}$$
$$R^o\lesssim R^i, \quad R^a\lesssim(\text{DoF})^{1/d}, \quad (R^a)^{-1}\lesssim(\text{DoF})^{-1/d}(\log\text{DoF})^{1/d},$$

for BQCE, P1-BQCF, and P1-BGFC $R^i\leq C_\Omega(R^a)^{1+2/d}$,

for P2-BQCF, P2-BGFC $R^i\leq C_\Omega(R^a)^{1+4/d}$.

In the assumption and in the rest of the paper, we write $A\lesssim B$ if there exists a constant $C$ such that $|A|\leq CB$, where $C$ is independent of $(\beta,\mathcal{T}_h)$, but may depend on the constants $\mathbf{C}$, or on any specified functions involved in the estimate (in particular, $C$ may depend on a solution $u^a$ and on derivatives of $V(Du^a)$ in some specified range, but not on test functions).

**Theorem 2.1** (Error estimates for blended methods). *For any given set of regularity constants $\mathbf{C}$ of $(\beta,\mathcal{T}_h)$, there exists $R_0^a>0$ such that, for all $(\beta,\mathcal{T}_h)$ satisfying Assumption 1, and in addition $R^a\geq R_0^a$, there exist solutions to BQCE (2.7), P1-BGFC (2.8), P2-BGFC (2.11), P1-(P2-)BQCF (2.13) methods, such that,*

$$\left\|\nabla u_h^{bqce}-\nabla\bar{u}^a\right\|_{L^2(\mathbb{R}^d)}\lesssim E^{\text{apx}}(u^a)+E^{\text{cb}}(u^a)+E^{\text{int}}(u^a)\lesssim(\text{DoF})^{1/2-2/d}, \tag{2.20}$$

$$\left\|\nabla u_h^{bgfc}-\nabla\bar{u}^a\right\|_{L^2(\mathbb{R}^d)}\lesssim E^{\text{apx}}(u^a)+E^{\text{cb}}(u^a)+E^{\text{int},2}(u^a)\lesssim(\text{DoF})^{-1/2-1/d}, \tag{2.21}$$

$$\left\|\nabla u_h^{bgfc,2} - \nabla \bar{u}^a\right\|_{L^2(\mathbb{R}^d)} \lesssim E^{apx,2}(u^a) + E^{cb}(u^a) + E^{int,2}(u^a) \lesssim (\mathrm{DoF})^{-1/2-2/d}, \quad (2.22)$$

$$\left\|\nabla u_h^{bqcf} - \nabla \bar{u}^a\right\|_{L^2(\mathbb{R}^d)} \lesssim E^{apx}(u^a) + E^{cb}(u^a) \lesssim (\mathrm{DoF})^{-1/2-1/d}, \quad (2.23)$$

$$\left\|\nabla u_h^{bqcf,2} - \nabla \bar{u}^a\right\|_{L^2(\mathbb{R}^d)} \lesssim E^{apx,2}(u^a) + E^{cb}(u^a) \lesssim (\mathrm{DoF})^{-1/2-2/d}. \quad (2.24)$$

We will sketch the proof for the first part of the inequalities in Theorem 2.1 in Section 4, and discuss the second part, namely, the error contributions in terms of DoF in the following paragraphs. Namely, the error contributions $E^{apx}$, $E^{cb}$, $E^{int}$ et. al. can be optimized with respect to DoF in terms of mesh parameters $(R^a, R^b, R^i, R^o, \beta, h(x))$.

We have the regularity results for the point defect configuration [10, Theorem 1], namely, the displacement field has the following decay away from the defect site:

$$\begin{aligned} |\nabla^j \tilde{u}^a(x)| &\lesssim |x|^{1-d-j} \quad \text{for } j = 0, \ldots, 3, \text{ and} \\ |\nabla^j \bar{u}^a(x)| &\lesssim |x|^{1-d-j} \quad \text{for } j = 0, 1. \end{aligned} \quad (2.25)$$

By (2.25) and Assumption 1, when $R^a$ is large, we have

$$\|\nabla \bar{u}^a\|_{L^2(\Omega^{ext})} \sim (R^i)^{-d/2}, \quad \|h\nabla^2 \tilde{u}^a\|_{L^2(\Omega^c)} \sim (R^a)^{-d/2-1},$$
$$\|h^2 \nabla^3 \tilde{u}^a\|_{L^2(\Omega^c)} \sim (R^a)^{-d/2-2}, \quad \|\nabla^3 \tilde{u}^a\|_{L^2(\Omega^c)} \sim (R^a)^{-d/2-2}, \quad (2.26)$$
$$\|\nabla^2 \tilde{u}^a\|_{L^4(\Omega^c)}^2 \sim (R^a)^{-3d/2-2}, \quad \|\nabla^2 \beta\|_{L^2(\Omega^b)} \|\nabla \bar{u}^a\|_{L^2(\Omega^b)} \sim (R^a)^{d/2-2}(R^a)^{-d} = (R^a)^{-d/2-2}.$$

We are in the position to estimate the errors of blended methods separately.

**BQCE method:** Picking the leading order terms from the error estimate (2.20), we have

$$\left\|\nabla u_h^{bqce} - \nabla \bar{u}^a\right\|_{L^2(\mathbb{R}^d)} \lesssim \|\nabla^2 \beta\|_{L^2(\Omega^b)} + \|\beta h \nabla^2 \tilde{u}^a\|_{L^2(\Omega_h)} + \|\nabla \bar{u}^a\|_{L^2(\Omega^{ext})} + \text{h.o.t.} \quad (2.27)$$

The term $\|\nabla^2 \beta\|_{L^2}$ in the coupling error $E^{int}$ is due to the smeared ghost forces, and an optimal choice of $\beta$ yields $\|\nabla^2 \beta\|_{L^2} \sim (R^a)^{d/2-2}$. The term $\|\beta h \nabla^2 \tilde{u}^a\|_{L^2(\Omega_h)}$ measures the finite element approximation error, while the term $\|\nabla \bar{u}^a\|_{L^2(\Omega^{ext})}$ measures the truncation error for a finite computational domain. Both of them are of the order $(R^a)^{-d/2-1}$ by (2.26) and the relation $R^i \sim (R^a)^{1+2/d}$. Therefore, the interface error is dominant:

$$\left\|\nabla u_h^{bqce} - \nabla \bar{u}^a\right\|_{L^2(\mathbb{R}^d)} \lesssim (R^a)^{d/2-2} \lesssim (\mathrm{DoF})^{1/2-2/d}. \quad (2.28)$$

**P1-BGFC method:** The analysis in Section 4 will reveal that the dominant coupling error contribution $\nabla^2 \beta$ comes from terms like $|\partial W(\nabla \tilde{u}) \nabla^2 \beta|$ or $|\partial V_\ell(Du) \nabla^2 \beta|$, where $W$ and $V$ can be replaced their renormalized counterparts (as in BQCE), or second renormalized counterparts (as in BGFC).

For the BQCE method, we can only bound $|\partial W'(\nabla \tilde{u})|$ or $|\partial V'_\ell(Du)|$ by their $L^\infty$ bounds. For the BGFC method, the following property will help us improve the estimate

$$\partial V''_\ell(0) = 0, \quad \partial W''(0) = 0. \tag{2.29}$$

We have

$$|\partial V_\ell(Du)\nabla^2\beta| = |\partial(V''_\ell(Du) - V''_\ell(0))\nabla^2\beta| \lesssim |\nabla \tilde{u}||\nabla^2\beta|, \tag{2.30}$$

namely, the coupling error is now combined with $\nabla \bar{u}$, which decays with the order $(R^a)^{-d/2-2}$, when $\bar{u}$ is replaced with $\bar{u}^a$ and as $R^a \to \infty$ in (2.21). For the P1-BGFC method, the coupling error is no longer dominant, instead the term $\|\beta h \nabla^2 \tilde{u}^a\|_{L^2(\Omega_h)}$ in the best-approximation error $E^{\text{apx}}$ becomes the leading contribution, which is of the order $(R^a)^{-d/2-1}$. Taking $R^i \sim (R^a)^{1+2/d}$ to balance the approximation error with the truncation error term $\|\nabla \bar{u}^a\|_{L^2(\Omega^{\text{ext}})}$, we have

$$\left\|\nabla u_h^{\text{bgfc}} - \nabla \bar{u}^a\right\|_{L^2(\mathbb{R}^d)}$$
$$\lesssim \|\nabla^2 \beta\|_{L^2(\Omega^b)} \|\nabla \tilde{u}\|_{L^2(\Omega^b)} + \|\beta h \nabla^2 \tilde{u}^a\|_{L^2(\Omega_h)} + \|\nabla \bar{u}^a\|_{L^2(\Omega^{\text{ext}})} + \text{h.o.t.}$$
$$\lesssim (R^a)^{-d/2-2} + (R^a)^{-d/2-1} + (R^i)^{-d/2} \lesssim (R^a)^{-d/2-1} \lesssim (\text{DoF})^{-1/2-1/d}.$$

**P2-BGFC method:** P1-BGFC is still sub-optimal in the sense that the term $\|\beta h \nabla^2 \tilde{u}^a\|_{L^2(\Omega_h)}$ in the best-approximation error dominates. This can be improved by using P2 approximation space $\mathscr{U}_h^{(2)}$ defined in (2.10). Now the leading terms in (2.22) are

$$\left\|\nabla u_h^{\text{bgfc},2} - \nabla \bar{u}^a\right\|_{L^2(\mathbb{R}^d)} \tag{2.31}$$
$$\lesssim \|\nabla^2 \beta\|_{L^2(\Omega^b)} \|\nabla \tilde{u}\|_{L^2(\Omega^b)} + \|h^2 \nabla^3 \tilde{u}^a\|_{L^2(\Omega_h)} + \|\nabla^3 \tilde{u}^a\|_{L^2(\Omega_h)} + \|\nabla \bar{u}^a\|_{L^2(\Omega^{\text{ext}})} + \text{h.o.t.}$$
$$\lesssim (R^a)^{-d/2-2} + (R^a)^{-d/2-2} + (R^a)^{-d/2-2} + (R^i)^{-d/2} \lesssim (R^a)^{-d/2-2} \lesssim (\text{DoF})^{-1/2-2/d}.$$

where the contributions from $E^{\text{apx}}$, $E^{\text{cb}}$, and $E^{\text{int}}$ are all balanced by taking $R^i \sim (R^a)^{1+4/d}$, and in that sense, this error estimate cannot be further improved if we use the Cauchy-Born model as the continuum model.

**BQCF method:** As a force based method, there is no coupling error in the BQCF error estimate (2.23). Therefore, it has the following convergence rate, with the same parameters as the BQCE method,

$$\left\|\nabla u_h^{\text{bqcf}} - \nabla \bar{u}^a\right\|_{L^2(\mathbb{R}^d)} \lesssim (R^a)^{-d/2-1} \lesssim (\text{DoF})^{-1/2-1/d}. \tag{2.32}$$

For P2-BQCF method, the approximation error can be improved by taking $\mathscr{U}_h^{(2)}$ in (2.13), similar to P2-BGFC. And the convergence rate is

$$\left\|\nabla u_h^{\text{bqcf},2} - \nabla \bar{u}^a\right\|_{L^2(\mathbb{R}^d)} \lesssim (R^a)^{-d/2-2} \lesssim (\text{DoF})^{-1/2-2/d}. \tag{2.33}$$

**Remark 2.1.** The energy error estimate for BQCE and P1-BGFC can be seen in [15, Theorem 3] [27, Theorem 5.1, Remark 4]. In general, we expect such a convergence rate,

$$\left|\mathcal{E}^{\text{a}}(u^{\text{a}}) - \mathcal{E}_h^{\beta}\left(u_h^{\beta}\right)\right| \lesssim \left\|\nabla u_h^{\beta} - \nabla \bar{u}^{\text{a}}\right\|_{L^2(\mathbb{R}^d)}^2. \tag{2.34}$$

P2-BGFC has the same convergence rate as P1-BGFC for the energy error.

**Remark 2.2.** Notice that when $u_0 = 0$, and $\Lambda = \Lambda^{\text{hom}}$, we have

$$\mathcal{E}^{\text{bgfc}}(u_h) = \mathcal{E}^{\text{bqce}}(u_h) - \left\langle \delta\mathcal{E}^{\text{bqce}}(0), u_h \right\rangle = \mathcal{E}^{\text{bqce}}(u_h) - \left\langle \delta\mathcal{E}^{\text{bqce}}(0) - \mathcal{F}_h^{\text{bqcf}}(0), u_h \right\rangle, \tag{2.35}$$

where $\mathcal{F}^{\text{bqcf}}$ is the BQCF operator defined in (2.12), and $\mathcal{F}^{\text{bqcf}}(0) = 0$ since it has no ghost forces. The second renormalization in BGFC is therefore equivalent to the dead load ghost-force correction scheme of Shenoy et. al. [32], applied for a blended coupling formulation and in the reference configuration.

BGFC scheme can be generalized as the following predictor-corrector formulation:

$$\mathcal{E}^{\text{bgfc}}(u_h) := \mathcal{E}^{\text{bqce}}(u_h) - \left\langle \delta\mathcal{E}^{\text{bqce}}(u_0) - \mathcal{F}_h^{\text{bqcf}}(u_0), u_h - u_0 \right\rangle, \tag{2.36}$$

where $u_0$ is a suitable reference configuration, or "predictor", that can be cheaply obtained. This kind of view gives use more flexibility in applications for cracks and dislocations.

## 3 Numerical Experiments

### 3.1 Interaction Potential

Many practical site potentials has the form $V(Du)$, for example, the generic pair functional form [34]. That includes the widely used embedded atom model (EAM) [4] and Finnis-Sinclair model [11]. Namely, the potential is a function of the distances between atoms within the interaction range and with no angular dependence. In our numerical implementation, the site potential is given by a toy EAM model (3.1), for which $V_\ell$ is of the form

$$V_\ell(Du) = \frac{1}{2} \sum_{\rho \in \mathcal{R}_\ell} \phi(|D_\rho u(\ell) + \rho|) + F\left(\sum_{\rho \in \mathcal{R}_\ell} \psi(|D_\rho u(\ell) + \rho|)\right), \tag{3.1}$$

$$\text{with} \quad \phi(r) = [e^{-2a(r-1)} - 2e^{-a(r-1)}], \quad \psi(r) = e^{-br},$$

$$F(\tilde{\rho}) = c\left[(\tilde{\rho} - \tilde{\rho}_0)^2 + (\tilde{\rho} - \tilde{\rho}_0)^4\right],$$

and with parameters $a = 4.4$, $b = 3$, $c = 5$, $\tilde{\rho}_0 = 12e^{-b}$. We consider next nearest neighbor (in hopping distance) interactions.

## 3.2 Mesh Generation

In this section, we provide more details for the mesh generation of the blended model. Let us assume that we have an FCC lattice with a point vacancy at the origin in three dimensions. The construction can be extended to other lattices, and also to separated vacancies and microcrack. The nodes $\mathcal{X}_h$ of the mesh $\mathcal{T}_h$ are layers of points with each layer forming the shell of an octahedron, which in turn forms a graded sequence of tetrahedra, see Figure 1(a). We assume the symmetry of $\mathcal{X}_h$ with respect to axes, therefore we only need to consider the positive octant $\mathbb{R}^3_+$ while other seven parts can be constructed by symmetry. Let $C(n,r)$ with $n \in \mathbb{N}$, $r \in \mathbb{R}_+$, be a closed shell of nodes defined by

$$C(n,r) = \frac{r}{n} \left\{ (x,y,z) \in \mathbb{Z}^3 \,|\, |x|+|y|+|z|=n \right\}, \tag{3.2}$$

where $n+1 = \#\{C(n,r) \cap \{(x,y,z)|x \geq 0, y \geq 0, z = 0\}\}$ is the number of nodes in the intersections of $C(n,r)$ with the positive $x-y$ plane (by symmetry, also for the positive $y-z$ and $x-z$ planes), $r$ is the distance between the positive $x-$intercept (by symmetry, also for the positive $y-$ and $z-$ intercept) of $C(n,r)$ and the origin. If $n=r$, we use the shorthand notation $C(r)$ for $C(n,r)$.

We focus on the tetrahedra mesh generation between two neighboring shells of nodes inside the positive octant $\mathbb{R}^3_+$. Assume $C^{(1)} = C(n^{(1)}, r^{(1)})$ is the inner layer of nodes, and $C^{(2)} = C(n^{(2)}, r^{(2)})$ is the outer layer of nodes, with $r^{(2)} > r^{(1)} > 0$. See Figure 3 for an illustration.

For all possible neighboring layers $C^{(1)}$ and $C^{(2)}$ satisfying the constraint $|n^{(1)} - n^{(2)}| \leq 2$, we can introduce a structured partition in the shell between $C^{(1)}$ and $C^{(2)}$, see Section 6.4 for details. The construction can be done in all octants consistently, and results in a closed shell of tetrahedra (also consistent with neighboring shell of tetrahedra).

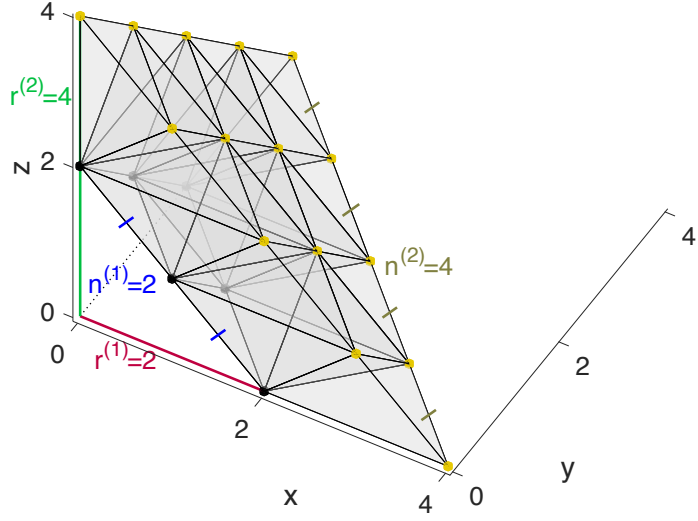

Figure 3: Partition in the atomistic region of an FCC lattice. We show the shell between two neighboring layers in the [111] direction, and in the positive octant $\mathbb{R}^3_+$.

In the atomistic region $\Omega^{\rm a}$ (including the blending region), the atomistic shells are $C(2k,2k)$, with $2k \leq r^{\rm a}$ and $k \in \mathbb{Z}_+$. $r^{\rm a}$ is the outer radius of $\Omega^{\rm a}$ which we assume to be an even integer, which is more convenient to use in the numerical implementation compared to $R^{\rm a}$, the inner radius of the set $\mathrm{supp}(\beta=0)$. We have the following relation:

$$R^{\rm a} = 2\left\lceil \left(r^{\rm a}+2-2R^{\rm cut}+R^{\rm def}\right)/4 \right\rceil / \sqrt{d} - 2R^{\rm cut} - \sqrt{d}.$$

In the continuum region $\Omega^{\rm c}$, distance between neighbor layers are determined according to the mesh size function $h(x)$ defined in (2.19). We start from several atomistic shells $C(r_i)$, with $n_1 = r_1 = r^{\rm a}$, $n_2 = r_2 = r^{\rm a}+2$, $n_3 = r_3 = r^{\rm a}+4$, and define a sequence $r_i$, such that $r_{i+1} - r_i \sim r_i^\chi$, where $\chi$ is the exponent in $h(x)$. We generate the shell of nodes $C(n_i, r_i)$ for $i > 3$, with $n_i$ defined by

$$n_i = \underset{n \in \mathbb{Z}_+, |n-n_{i-1}| \leq 2}{\mathrm{argmin}} |n(r_i - r_{i-1}) - r_i|.$$

We generate $M$ shells of nodes in the continuum region such that $r_{M-1} < r^{\rm c} \leq r_M$, where $r^{\rm c} := R^{\rm a} + (R^{\rm a})^{5/3}$ for P1-BGFC, BQCE and P1-BQCF, and $r^{\rm c} := R^{\rm a} + (R^{\rm a})^{7/3}$ for P2-BGFC and P2-BQCF, also see Assumption 1.

### 3.3 Single Vacancy Example

We first consider the case with one single vacancy at the origin. We use a reference solution with atomistic radius $r^{\rm a} = 60$ by P1-BGFC. We numerically solve P1-BGFC and

P1-BQCF solutions with $r^a = 8, 12, \cdots, 56$, BQCE solutions with $r^a = 8, 12, \cdots, 36$, and P2-BGFC and P2-BQCF solutions with $r^a = 8, 12, \cdots, 24$. We apply an uniform deformation at the boundary of the computational domain: $y|_{\partial \Omega_h} = y_\mathsf{F}$, where,

$$\mathsf{F} = \begin{pmatrix} 1 & 0.01 & 0.02 \\ 0 & 1 & 0.015 \\ 0 & 0 & 1 \end{pmatrix}. \tag{3.3}$$

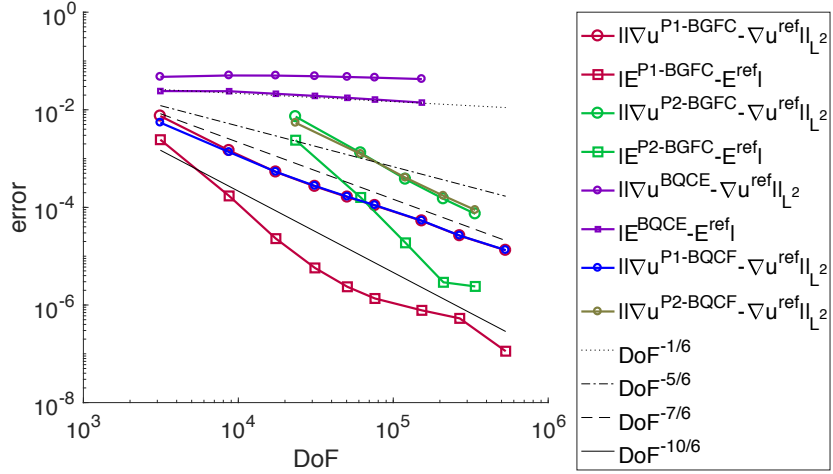

(a) Error vs. number of DoF.

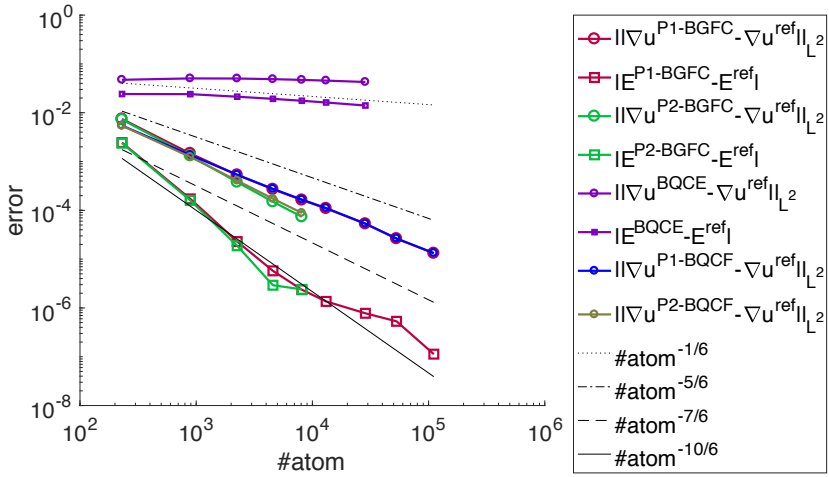

(b) Error vs. number of atoms.

Figure 4: Convergence error (P1-BGFC, P2-BGFC, BQCE, P1-BQCF, and P2-BQCF) for point vacancy. We use the P1-BGFC reference solution with $r^a = 60$.

We use the Newton-Raphson method to solve the blended coupling solutions in (2.8) etc., and we stop the iteration when the $\ell^\infty$ norm of the increment error is less than $10^{-14}$. We use a Linux cluster with Intel Xeon E5 x86-64 CPU with 256 cores and 2TB memory. The program is written in MATLAB.

The $H^1$ semi-norm errors for BGFC (P1- and P2-) together with BQCE and BQCF (P1- and P2-) are plotted with respect to the DoF in Figure 4(a), and with respect to #atom (DoF in the atomistic region) in Figure 4(b). In Figure 4(a), the P1-BGFC displacement error decays faster than the theoretical prediction $\mathcal{O}\left(\text{DoF}^{-5/6}\right)$, and P1-BQCF has a very similar curve. The P2-BGFC displacement error also decays faster than the theoretical prediction $\mathcal{O}\left(\text{DoF}^{-7/6}\right)$, but with a larger constant compared to P1-BGFC. The curves of P1- and P2- BGFC do not cross over until the maximum DoF of about $5\times 10^5$ (using about 1.8 TB memory) is reached. Again, P2-BQCF has a very similar curve compared to P2-BGFC. Following the similar analysis [5] in 2D, the improvement of P2-BGFC can be seen when plotting error with respect to #atom in Figure 4(b).

We numerically show the decay of displacement errors vs. the distance to the defect core for solutions with fixed $r^a$. For example, $r^a = 28$ for the P1-BGFC and P2-BGFC solutions in Figure 5. We can see that P2-BGFC significantly reduces the error contribution in the continuum region. And for the solutions with large $r^a$, we can observe that the main error is in the continuum region, so P2-BGFC causes more error reduction than P1-BGFC does. We also observe that all numerical solutions observe the decay property (2.25) for $d=3$ and $j=0$.

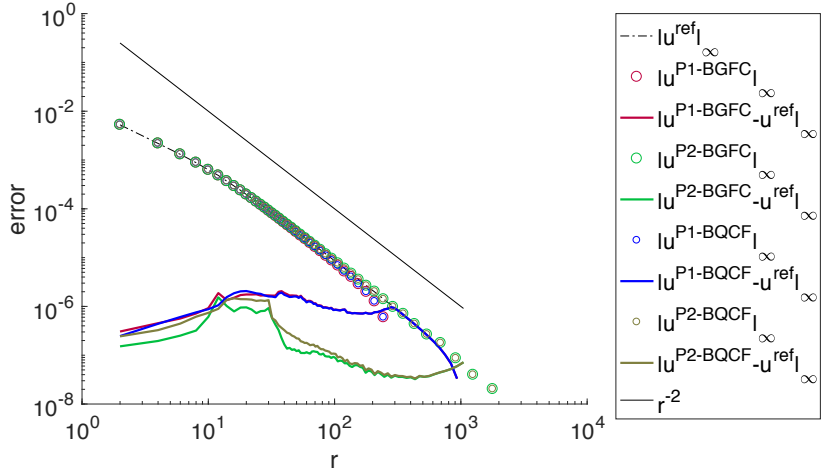

Figure 5: Decay of displacement errors in $\ell^\infty$ norm with respect to the distance to the defect core, for approximate solutions with $r^a = 28$. The distance $r$ is taken with respect to the adjacent shell of nodes $C(n,r)$. The blending region for the approximate solution is $12 \leq r \leq 26$, which has larger error. In the legend, $u^{\text{ref}}$ represents the reference solution on the reference mesh, $u^\beta$ represents the approximate solution on a coarser mesh or the interpolated approximate solution on the reference mesh.

### 3.4 Separated Vacancies and Microcrack Example

We also consider other defect types, such as two separated point vacancies at $(x,0,0)$, $x = \pm 4$; and a micro-crack consisting of five vacancies in a row at $(x,0,0)$, $x = 0, \pm 2, \pm 4$ in the reference configuration. In these examples, we use a P2-BGFC reference solution with $r^a = 28$, and the approximation solutions are P1-BGFC solutions with $r^a = 12, 16, \cdots, 56$ and P2-BGFC solutions with $r^a = 12, 16, 20, 24$. We show the $H^1$ semi-norm errors in Figure 6. Comparing these two examples and the single-vacancy example, we find that the convergence curve has a larger constant with respect to the size of the defect core.

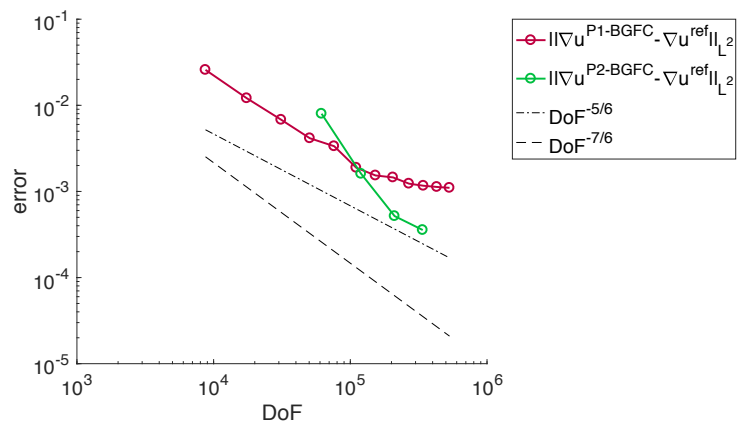

(a) Separated-vacancy case.

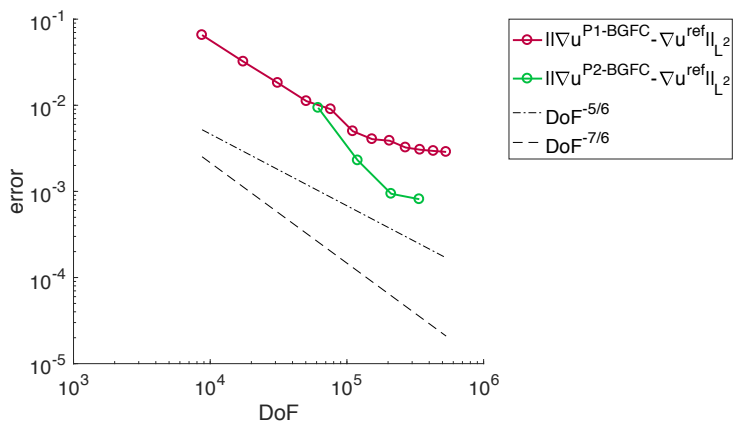

(b) Micro-crack case.

Figure 6: $H^1$ semi-norm error for P1-BGFC and P2-BGFC, with a P2-BGFC reference solution with $r^a = 28$.

## 4 Analysis

In this section, we briefly sketch the error analysis of blending type methods such as BQCE, BQCF, and BGFC under a unified analytical framework. We will emphasize the difference and connection of those methods. More technical details can be found in [15, 27].

### 4.1 Framework

We adopt the analytical framework for a/c coupling methods in [10, 18], which is analogous to that of finite element methods for regular nonlinear PDE, employing quasi-best approximation, consistency and stability.

Let $\mathscr{G}_h := \delta\mathscr{E}_h^\beta + \delta\mathscr{P}^{\text{def}}$ for the blended coupling problem be the force operator of BQCE or BGFC scheme, or $\mathscr{G}_h := \mathscr{F}_h^\beta + \delta\mathscr{P}$ for the BQCF scheme. Let $\Pi_h : \mathscr{U}^{1,2} \to \mathscr{U}_h$ be the quasi-best approximation operator which maps the displacement from the atomistic lattice to the coarse-grained mesh, see [15, Section 4.2.4], such that, $\|\nabla\Pi_h u - \nabla\bar{u}\|_{L^2} \lesssim E^{\text{apx}}$, with $E^{\text{apx}}$ defined in (2.14).

Assume $u^{\text{a}}$ is strongly stable with stability constant $\gamma_{\text{a}} > 0$ in (2.5), we shall require that $\mathscr{G}_h$ is consistent for some small consistency error $\eta > 0$ that depends on $u^{\text{a}}$, $\mathcal{T}_h$, and $\beta$,

$$\forall v_h \in \mathscr{U}_h, \quad \langle \mathscr{G}_h(\Pi_h u^{\text{a}}), v_h \rangle \leq \eta \|\nabla v_h\|_{L^2(\mathbb{R}^d)}, \tag{4.1}$$

and stable,

$$\forall v_h \in \mathscr{U}_h, \quad \langle \delta\mathscr{G}_h(\Pi_h u^{\text{a}}) v_h, v_h \rangle \geq c_0 \|\nabla v_h\|_{L^2(\mathbb{R}^d)}^2. \tag{4.2}$$

In general, we expect that with some function $\omega(R^{\text{a}}) \to 0$ as $R^{\text{a}} \to \infty$, it holds true that $c_0 \geq \gamma_{\text{a}} - \omega(R^{\text{a}})$ [15, 27].

We then employ the Inverse Function Theorem [18] to prove that, if $\eta/c_0$ is sufficiently small (adding some technical assumptions), there exists $w_h \in \mathscr{U}_h$ such that

$$\|\nabla w_h\|_{L^2(\mathbb{R}^d)} \leq 2\eta/c_0, \quad \mathscr{G}_h(\Pi_h u^{\text{a}} + w_h) = 0.$$

Thus, we can construct a blended coupling solution $u_h^\beta := \Pi_h u^{\text{a}} + w_h$ satisfying $\mathscr{G}_h\left(u_h^\beta\right) = 0$, and,

$$\begin{aligned}\left\|\nabla u_h^\beta - \nabla\bar{u}^{\text{a}}\right\|_{L^2(\mathbb{R}^d)} &\leq \left\|\nabla u_h^\beta - \nabla\Pi_h u^{\text{a}}\right\|_{L^2(\mathbb{R}^d)} + \|\nabla\Pi_h u^{\text{a}} - \nabla\bar{u}^{\text{a}}\|_{L^2(\mathbb{R}^d)} \\ &\leq \|\nabla w_h\|_{L^2(\mathbb{R}^d)} + \|\nabla\Pi_h u^{\text{a}} - \nabla\bar{u}^{\text{a}}\|_{L^2(\mathbb{R}^d)} \\ &\leq 2\frac{\eta}{c_0} + \|\nabla\Pi_h u^{\text{a}} - \nabla\bar{u}^{\text{a}}\|_{L^2(\mathbb{R}^d)}.\end{aligned} \tag{4.3}$$

The first term in (4.3), $2\eta/c_0$ is sufficiently small when $R^{\text{a}}$ is sufficiently large and $(\beta, \mathcal{T}_h)$ satisfies Assumption 1, and the second term can be bounded by the quasi-best approximation error $E^{\text{apx}}$.

## 4.2 Sketch of the proof

We focus on the analysis of BQCE and BGFC schemes. The analysis for P1-BQCF can be found in [15], while the P2-BQCF convergence is the straightforward extension of P1-BQCF results in addition to the approximation property of $\mathscr{U}_h^{(2)}$.

We have a unified representation of the blended energy $\mathscr{E}^\beta$ for BQCE and BGFC. For $u_h \in \mathscr{U}_h$, we define

$$\mathscr{E}_h^\beta(u_h) := \sum_{\ell \in \Lambda^{\text{hom}} \cap \Omega_h} (1-\beta(\ell))V_\ell(Du) + \int_{\Omega_h} Q_h[\beta W_\ell(\nabla u_h)]. \tag{4.4}$$

We obtain the BQCE energy (2.6) by replacing $V$, $W$ with $V'$, $W'$, and obtain the P1-BGFC energy (2.8) by replacing $V$, $W$ with $V''$, $W''$. For P2-BGFC, we can replace $\mathscr{U}_h$ by $\mathscr{U}_h^{(2)}$ in the P1-BGFC energy.

We now derive the atomistic stress by the weak form. For $\forall u, v \in \mathscr{U}^{1,2}$,

$$\langle \delta\mathscr{E}^{\text{a}}(u), v^* \rangle = \sum_{\ell \in \Lambda^{\text{hom}}} \sum_{\rho \in \mathcal{R}} V_{\ell,\rho} \cdot D_\rho v^*(\ell) \tag{4.5}$$

$$= \sum_{\ell \in \Lambda^{\text{hom}}} \sum_{\rho \in \mathcal{R}} V_{\ell,\rho} \cdot \int_{\mathbb{R}^d} \omega_\rho(\ell-x) \nabla_\rho \bar{v} \, dx \tag{4.6}$$

$$= \int_{\mathbb{R}^d} \left\{ \sum_{\ell \in \Lambda^{\text{hom}}} \sum_{\rho \in \mathcal{R}} [V_{\ell,\rho} \otimes \rho] \omega_\rho(\ell-x) \right\} : \nabla \bar{v}(\ell), \tag{4.7}$$

notice that we take test function $v^*$ instead of $v$, which is defined in the Section 6.1, and use the property of $\omega_\rho(\ell-x)$ in (6.1). We can define the corresponding atomistic stress as $S^{\text{a}}(y;x) := \sum_{\ell \in \Lambda^{\text{hom}}} \sum_{\rho \in \mathcal{R}} [V_{\ell,\rho} \otimes \rho] \omega_\rho(\ell-x)$. The region $\nu_x$ is the affected neighbourhood of some $x \in \mathbb{R}^d$ that contains all lattice points $\ell \in \Lambda^{\text{hom}}$ involve in the atomistic stress

$$\forall x \in \mathbb{R}^d, \quad \nu_x := B_{2R^{\text{cut}}+\sqrt{d}}(x). \tag{4.8}$$

Given $\forall u_h, v_h \in \mathscr{U}_h$, we define $u := u_h|_{\Lambda^{\text{hom}}} \in \mathscr{U}^{1,2}$, and $v := \Pi_h' v_h \in \mathscr{U}^{\text{c}}$ through the dual approximation operator defined in Section 6.3. We note that the first and third terms in (4.9) cancel out due to the fact that $u(\ell) = u_h(\ell)$, $v^*(\ell) = v_h(\ell)$, $\forall \ell \in \Lambda^{\text{a}}$.

$$\left\langle \delta\mathscr{E}_h^\beta(u_h), v_h \right\rangle - \left\langle \delta\mathscr{E}^{\text{a}}(u), v^* \right\rangle$$

$$= \sum_{\ell \in \Lambda^{\text{hom}}} (1-\beta(\ell)) \langle \delta V(Du_h(\ell)), Dv_h(\ell) \rangle + \int_{\mathbb{R}^d} Q_h[\beta(x) \partial W(\nabla u_h) : \nabla v_h] \, dx \tag{4.9}$$

$$- \sum_{\ell \in \Lambda^{\text{hom}}} (1-\beta(\ell)) \langle \delta V(Du(\ell)), Dv^*(\ell) \rangle - \int_{\mathbb{R}^d} \sum_{\ell \in \Lambda^{\text{hom}}} \beta(\ell) \sum_{\rho \in \mathcal{R}} \omega_\rho(\ell-x) V_{,\rho}(Du(\ell)) \otimes \rho : \nabla \bar{v} \, dx$$

$$= \int_{\mathbb{R}^d} Q_h((\beta(x)\partial W(\nabla u_h)) : \nabla v_h) \, dx - \int_{\mathbb{R}^d} \sum_{\ell \in \Lambda^{\text{hom}}} \beta(\ell) \sum_{\rho \in \mathcal{R}} \omega_\rho(\ell-x) V_{,\rho}(Du(\ell)) \otimes \rho : \nabla \bar{v} \, dx$$

$$= T_1 + T_2 + T_3 + T_4,$$

where

$$T_1 := \int_{\mathbb{R}^d} Q_h[(\beta(x)(\partial W(\nabla u_h) - \partial W(\nabla \tilde{u}))) : \nabla v_h] \, \mathrm{d}x, \tag{4.10}$$

$$T_2 := \int_{\mathbb{R}^d} (Q_h - 1)(\beta(x)\partial W(\nabla \tilde{u}) : \nabla v_h(x)) \, \mathrm{d}x, \tag{4.11}$$

$$T_3 := \int_{\mathbb{R}^d} (\beta(x)\partial W(\nabla \tilde{u})) : (\nabla v_h(x) - \nabla \bar{v}) \, \mathrm{d}x, \tag{4.12}$$

$$T_4 := \int_{\mathbb{R}^d} \left( \beta(x)\partial W(\nabla \tilde{u}) - \sum_{\ell \in \Lambda^{\mathrm{hom}}} \beta(\ell) \sum_{\rho \in \mathcal{R}} \omega_\rho(\ell - x) V_{,\rho}(Du(\ell)) \otimes \rho \right) : \nabla \bar{v} \, \mathrm{d}x. \tag{4.13}$$

In the last term $T_4$, we define the stress error for $\tilde{u} \in C^\infty(\mathbb{R}^d; \mathbb{R}^d)$ as

$$\mathsf{R}^\beta(\tilde{u}; x) := \beta(x)\partial W(\nabla \tilde{u}) - \sum_{\ell \in \Lambda^{\mathrm{hom}}} \beta(\ell) \sum_{\rho \in \mathcal{R}} \omega_\rho(\ell - x) V_{,\rho}(Du(\ell)) \otimes \rho, \tag{4.14}$$

where $\partial W(\nabla \tilde{u})$ can be seen as the stress of the Cauchy-Born model, and the second term in $\mathsf{R}^\beta(\tilde{u};x)$ is the atomistic stress when $\beta(\ell) = 1$. Therefore, $T_4 = \int_{\mathbb{R}^d} \mathsf{R}^\beta(\tilde{u};x) : \nabla \bar{v} \, \mathrm{d}x$.

$T_1$ is bounded by the interpolation error of $\tilde{u}$ and $\bar{u}$, hence by the best approximation error $E^{\mathrm{apx}}$.

$$T_1 \lesssim \left( \|\nabla u_h - \nabla \tilde{u}\|_{L^2(\Omega^c)} + \|h^2 \nabla^3 \tilde{u}\|_{L^2(\Omega^c)} \right) \|\nabla v_h\|_{L^2(\Omega^c)}$$
$$\lesssim \left( \|\nabla \bar{u}\|_{L^2(\Omega^{\mathrm{ext}})} + \|h \nabla^2 \tilde{u}\|_{L^2(\Omega^c)} + \|h^2 \nabla^3 \tilde{u}\|_{L^2(\Omega^c)} \right) \|\nabla v_h\|_{L^2(\Omega^c)}. \tag{4.15}$$

If $\mathscr{U}_h$ is replaced by $\mathscr{U}_h^{(2)}$ as in P2-BGFC, we have

$$T_1 \lesssim \left( \|\nabla \bar{u}\|_{L^2(\Omega^{\mathrm{ext}})} + \|h^2 \nabla^3 \tilde{u}\|_{L^2(\Omega^c)} \right) \|\nabla v_h\|_{L^2(\Omega^c)}. \tag{4.16}$$

$T_2$ is the quadrature error of the midpoint rule:

$$T_2 \lesssim \|h^2 \nabla^2 (\beta \partial W(\nabla \tilde{u}))\|_{L^2(\Omega_h)} \|\nabla v_h\|_{L^2(\Omega_h)}. \tag{4.17}$$

$T_3$ can be estimated by the integration by parts, and the fact that for $\ell \in \mathcal{N}_h \cap \Lambda^a$, $v_h(\ell) = v^*(\ell)$, by the definition of dual approximation operator $\Pi'_h$ in Section 6.3. Together with the fact that $u_0 = 0$ for the point defect case, we have

$$T_3 \lesssim \left( \|\nabla^2 (\beta \partial W(\nabla \tilde{u}))\|_{L^2(\Omega^c)} + \|h \nabla^2 \tilde{u}\|_{L^2(\Omega^c)} \right) \|\nabla v_h\|_{L^2(\Omega_h)}. \tag{4.18}$$

For P2-BGFC method, we have an improved estimate:

$$T_3 \lesssim \left( \|\nabla^2 (\beta \partial W(\nabla \tilde{u}))\|_{L^2} \right) \|\nabla v_h\|_{L^2}. \tag{4.19}$$

The chain rule leads to

$$|\nabla^2(\beta\partial W(\nabla\tilde{u}))|\lesssim|\beta\nabla^3\tilde{u}|+|\nabla\beta\nabla^2\tilde{u}|+|\partial W(\nabla\tilde{u})\nabla^2\beta|, \tag{4.20}$$

where $|\partial W(\nabla\tilde{u})\nabla^2\beta|$ is the leading order term.

We need a more detailed analysis for $T_4$:

$$\begin{aligned}
&\mathsf{R}^\beta(\tilde{u};x)\\
&=\beta(x)\partial W(\nabla\tilde{u})-\sum_{\ell\in\Lambda}\beta(\ell)\sum_{\rho\in\mathcal{R}}\omega_\rho(\ell-x)V_{,\rho}(D\tilde{u}(\ell))\otimes\rho\\
&=\beta(x)\partial W(\nabla\tilde{u})-\sum_{\ell\in\Lambda}(\beta(x)+\nabla\beta(x)\cdot(\ell-x)+\mathcal{O}(\delta_2))\sum_{\rho\in\mathcal{R}}\omega_\rho(\ell-x)\\
&\quad\cdot\left(V_{,\rho}(\nabla_\mathcal{R}\tilde{u}(x))+\sum_{\xi\in\mathcal{R}}^N V_{,\rho\xi}(\nabla_\mathcal{R}\tilde{u}(x))\cdot(D\tilde{u}(\ell)-\nabla_\mathcal{R}\tilde{u}(x))+\mathcal{O}(\epsilon_2^2)\right)\otimes\rho\\
&=\beta(x)\partial W(\nabla\tilde{u})-\sum_{\ell\in\Lambda}(\beta(x)+\nabla\beta(x)\cdot(\ell-x)+\mathcal{O}(\delta_2))\sum_{\rho\in\mathcal{R}}\omega_\rho(\ell-x)\\
&\quad\cdot\left(V_{,\rho}(\nabla_\mathcal{R}\tilde{u}(x))+\sum_{\xi\in\mathcal{R}}^N V_{,\rho\xi}(\nabla_\mathcal{R}\tilde{u}(x))\cdot\left(\nabla_\xi\nabla_{\ell-x}\tilde{u}(x)+\frac{1}{2}\nabla_\xi^2\tilde{u}(x)+\mathcal{O}(\epsilon_3)\right)+\mathcal{O}(\epsilon_2^2)\right)\otimes\rho\\
&=\beta(x)\left(\partial W(\nabla\tilde{u})-\sum_{\rho\in\mathcal{R}}V_{,\rho}(\nabla_\mathcal{R}\tilde{u}(x))\otimes\rho\right)+\frac{1}{2}\sum_{\rho\in\mathcal{R}}(V_{,\rho}(\nabla_\mathcal{R}\tilde{u}(x))\otimes\rho)\nabla\beta(x)\cdot\rho\\
&\quad-\beta(x)\sum_{\rho,\xi\in\mathcal{R}}^N V_{,\rho\xi}(\nabla_\mathcal{R}\tilde{u}(x))\cdot\left(-\frac{1}{2}\nabla_\rho\nabla_\xi\tilde{u}+\frac{1}{2}\nabla_\xi^2\tilde{u}\right)\otimes\rho\\
&\quad+\mathcal{O}\left(|\nabla\beta(x)||\nabla^2\tilde{u}(x)|+\epsilon_2^2+(|\partial V(\nabla_\mathcal{R}\tilde{u}(x))|+|\nabla^2\tilde{u}(x)|)\delta_2+\epsilon_3\right) \tag{4.21}\\
&=0+0+0+\mathcal{O}\left(|\nabla\beta(x)||\nabla^2\tilde{u}(x)|+\epsilon_2^2+(|\partial V(\nabla_\mathcal{R}\tilde{u}(x))|+|\nabla^2\tilde{u}(x)|)\delta_2+\epsilon_3\right),
\end{aligned}$$

where $\epsilon_j:=\|\nabla^j\tilde{u}\|_{L^\infty(\nu_x)}$ and $\delta_j:=\|\nabla^j\beta\|_{L^\infty(\nu_x)}$.

The first three terms in (4.21) vanish by the following arguments: the first one is due to the definition of $W$; the second and the third ones are due to the symmetry of set $\mathcal{R}$ and the symmetry of $V$: $V_{,-\rho}=-V_{,\rho},V_{,\rho\xi}=-V_{,(-\rho)\xi}=-V_{,\rho(-\xi)}$.

For BQCE method, we can replace $V$ by $V'$, and obtain that

$$\left\|\mathsf{R}^\beta(\tilde{u};\cdot)\right\|_{L^2}\lesssim\|\nabla\beta\|_{L^\infty(\Omega^b)}\|\nabla^2\tilde{u}\|_{L^2(\Omega^b)}+\|\nabla^2\beta\|_{L^2(\Omega^b)}+\|\nabla^2\tilde{u}\|_{L^4(\Omega^c)}^2+\|\nabla^3\tilde{u}\|_{L^2(\Omega^c)}.$$

For BGFC method, we can replace $V$ by $V''$, and the changes in $T_4$ are in the following term,

$$\begin{aligned}
&\int_{\Omega^b}\left(|\partial V''(\nabla\tilde{u}(x))|+|\nabla^2\tilde{u}(x)|\right)^2\|\nabla^2\beta\|_{L^\infty(\nu_x)}^2\,\mathrm{d}x\\
&=\int_{\Omega^b}\left(|\partial V''(\nabla\tilde{u}(x))-\partial V''(0)|+|\nabla^2\tilde{u}(x)|\right)^2\|\nabla^2\beta\|_{L^\infty(\nu_x)}^2\,\mathrm{d}x
\end{aligned}$$

$$\lesssim \int_{\Omega^{\mathrm{b}}} \big(|\nabla\tilde{u}(x)|+|\nabla^2\tilde{u}(x)|\big)^2 \big\|\nabla^2\beta\big\|_{L^2(v_x)}^2 \mathrm{d}x$$
$$\lesssim \big\|\nabla^2\beta\big\|_{L^2(\Omega^{\mathrm{b}})}^2 \Big(\big\|\nabla\tilde{u}\big\|_{L^2(\Omega^{\mathrm{b}})}^2 + \big\|\nabla^2\tilde{u}\big\|_{L^2(\Omega^{\mathrm{b}})}^2\Big).$$

where the leading contribution $\big\|\nabla^2\beta\big\|_{L^2(\Omega^{\mathrm{b}})}$ in BQCE is now replaced by $\big\|\nabla^2\beta\big\|_{L^2(\Omega^{\mathrm{b}})}\big\|\nabla\tilde{u}\big\|_{L^2(\Omega^{\mathrm{b}})}$ in BGFC.

## 5 Conclusion and Outlook

In this paper, we demonstrate the 3D implementation of the BGFC coupling method, which achieves the optimal convergence rate for multi-body interaction potentials, and general interfaces. P1-BGFC can reach the same order of convergence as the force based P1-BQCF method and the (theoretical) QNL type methods. P2-BGFC and P2-BQCF have an optimal convergence rate among all the methods using Cauchy-Born model as the continuum elastic model, while BGFC has the advantage of being an energy based method.

We also review the convergence theorem of the blended coupling methods such as BQCE and BQCF, together with BGFC. We emphasize the difference of those coupling methods, and point out that the BGFC method can reduce the coupling error with the second renormalization of the potential, and P2-BGFC can achieve the optimal order due to the fact that it balance out all the contributions from coupling, Cauchy-Born coarse graining, and best approximation errors.

We implement the BGFC methods for single vacancy, separated vacancies, and microcrack in the three dimensional FCC lattice and with finite range multi-body interaction potential. We observed the theoretical convergence rate numerically. P2-BGFC and P2-BQCF methods admit the best decay rate as predicted.

We plan to implement BGFC method for multi-lattices, and also more realistic defects such as void and dislocations in our future work. Another possibility is to find higher order predictors in the predictor-corrector formulation, which can results in methods with higher order convergence rates.

## 6 Appendix

We collect a list of technical tools in Section 6.1, Section 6.2, and Section 6.3. We also provide the detail of tetrahedra partition between two neighboring lays in one octant in Section 6.4.

### 6.1 $P_1$ interpolant and the convolution

The lattice $\Lambda^{\mathrm{hom}}$ naturally induces a simplicial micro-triangulation $\mathcal{T}^{\mathrm{a}}$, see for example Section 3.2 for a possible construction in 3D.

Let $\zeta \in W^{1,\infty}(\Lambda^{\rm hom};\mathbb{R})$ be the P1 nodal basis function associated with the origin; namely, $\zeta$ is piecewise linear with respect to $\mathcal{T}^{\rm a}$, and $\zeta(0)=1$ and $\zeta(\ell)=0$ for $\ell \neq 0$ and $\ell \in \Lambda^{\rm hom}$. The nodal interpolant of $v \in \mathscr{U}$ can be written as

$$\bar{v}(x) := \sum_{\ell \in \Lambda^{\rm hom}} v(\ell)\zeta(x-\ell).$$

We define $v^* := \bar{\zeta} * \bar{v}$ for $v \in \mathscr{U}$,

$$D_\rho v^*(\ell) = \int_{s=0}^{1} \nabla_\rho v^*(\ell+s\rho)ds = \int_{\mathbb{R}^d}\int_{s=0}^{1}\zeta(\xi+s\rho-x)\nabla_\rho \bar{v}(x)dsdx$$
$$= \int_{\mathbb{R}^d}\omega_\rho(\ell-x)\nabla_\rho \bar{v} dx \quad \text{where } \omega_\rho(x) := \int_{s=0}^{1}\zeta(x+s\rho)dx,$$

and we have the following properties

$$\sum_\ell \omega_\rho(\ell-x) = 1, \quad \sum_\ell (\ell-x)\omega_\rho(\ell-x) = -\frac{1}{2}\rho. \tag{6.1}$$

### 6.2 $C^{2,1}$ conforming multi-quintic interpolant

In the analysis, we need higher order interpolation to measure the regularity of atomistic displacements. It is possible to define a $C^{2,1}$-conforming multi-quintic interpolation as in [15]. For $v: \Lambda^{\rm hom} \to \mathbb{R}^m$ and $i=1,\cdots,d$, let $\mathsf{d}_i^0 v(\ell) := v(\ell)$; $\mathsf{d}_i^1 v(\ell) := \frac{1}{2}(v(\ell+e_i)-v(\ell-e_i))$ and $\mathsf{d}_i^2 v(\ell) := v(\ell+e_i)-2v(\ell)+(\ell-e_i)$. Lemma 2.1 in [15] states that, for each $\ell \in \Lambda^{\rm hom}$ there exists a unique multi-quintic function $\tilde{v}: \ell+[0,1]^d \to \mathbb{R}^m$ defined through the conditions

$$\partial_{x_1}^{\alpha_1}\cdots\partial_{x_d}^{\alpha_d}\tilde{v}(\ell') = \mathsf{d}_1^{\alpha_1}\cdots\mathsf{d}_d^{\alpha_d}v(\ell'), \qquad \forall \ell' \in \ell+\{0,1\}^d, \alpha \in \{0,1,2\}^d, |\alpha|_\infty \leq 2,$$

and moreover, $\tilde{v}$ satisfies $\left\|\nabla^j \tilde{v}\right\|_{L^p(\xi+(0,1)^d)} \leq C \left\|D_{\mathcal{R}}^j v\right\|_{\ell^p(\xi+\{-1,0,1,2\}^d)}$.

### 6.3 Dual approximation operator

The quasi-best approximation operator $\Pi_h$ in [15, Section 4.2.4] defines a mapping from the atomistic space $\mathscr{U}^{\rm c}$ to the coarse grained space $\mathscr{U}_h$. The dual approximation operator $\Pi_h': \mathscr{U}_h \to \mathscr{U}^{\rm c}$, s.t., $\forall v_h \in \mathscr{U}_h$, $(\Pi_h' v_h)^*(\ell) = v_h(\ell)$ can be defined such that, $\forall \ell \in \Lambda^{\rm a}$, and $\Pi_h' v_h(\ell) = \zeta * v_h(\ell), \forall \ell \in \mathbb{Z}^d \setminus \Lambda^{\rm a}$.

It is proved in [15] that the dual approximation operator $\Pi_h'$ is well-defined. Moreover, $\exists C > 0, \forall v_h \in \mathscr{U}_h$,

$$\left\|\nabla \left(\Pi_h' v_h\right)^*\right\|_{L^2} \leq \left\|\nabla \overline{\Pi_h' v_h}\right\|_{L^2} \leq C\|\nabla v_h\|_{L^2}, \quad \left\|v_h - \overline{\Pi_h' v_h}\right\|_{L^2} \leq C\|\nabla v_h\|_{L^2}. \tag{6.2}$$

## 6.4 Partition of a shell between two neighboring layers

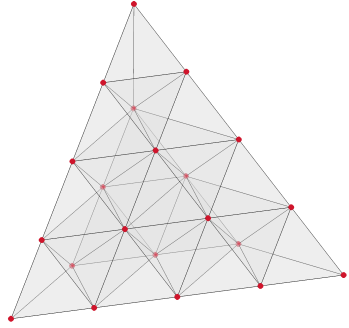

(a) $n^{(1)}=2$, $n^{(2)}=4$.

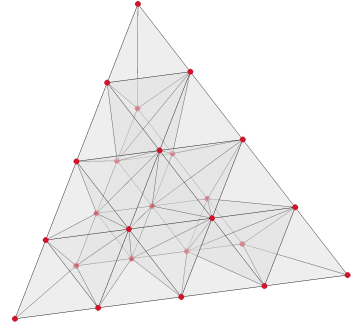

(b) $n^{(1)}=3$, $n^{(2)}=4$.

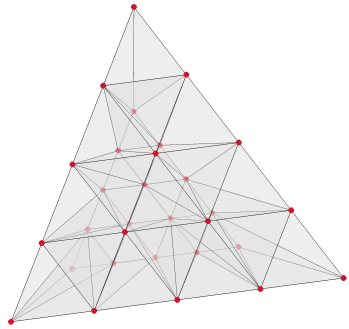

(c) $n^{(1)}=4$, $n^{(2)}=4$.

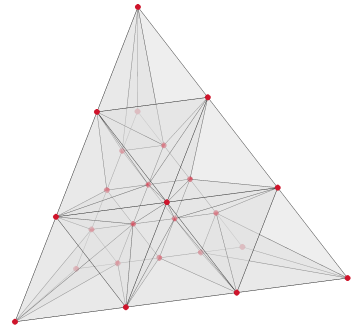

(d) $n^{(1)}=4$, $n^{(2)}=3$.

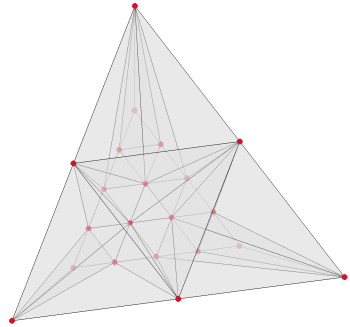

(e) $n^{(1)}=4$, $n^{(2)}=2$.

Figure 7: Examples of mesh partition in the continuum region.

Let $C^{(1)}$ and $C^{(2)}$ be two neighboring layers, with parameters $(n^{(1)},r^{(1)})$ and $(n^{(2)},r^{(2)})$ as defined in (3.2). We have the constraint $|n^{(1)}-n^{(2)}|\leq 2$. For each possible value of

$n^{(1)} - n^{(2)}$, we can introduce a structured partition of the shell between $C^{(1)}$ and $C^{(2)}$ as shown in Figure 7.

- if $n^{(2)} = n^{(1)}$, we construct $\left(n^{(1)}+1\right) n^{(1)}/2$ triangular prisms;

- if $n^{(2)} = n^{(1)}+1$, we construct $\left(n^{(1)}+2\right)\left(n^{(1)}+1\right)/2$ tetrahedra and $\left(n^{(1)}+1\right) n^{(1)}/2$ octahedra;

- if $n^{(2)} = n^{(1)}+2$, we construct $\left(n^{(1)}+1\right) n^{(1)}/2 + \left(n^{(1)}+2\right)\left(n^{(1)}+1\right)/2 + 3$ tetrahedra, $n^{(1)}\left(n^{(1)}-1\right)/2$ octahedra, and $3n^{(1)}$ pyramids;

- if $n^{(2)} = n^{(1)} - 1$ or $n^{(2)} = n^{(1)} - 2$, the construction is similar.

Those building blocks (prisms, octahedra, pyramids) can be further divided into tetrahedra.

## Acknowledgments

26



\end{document}